\documentclass[12pt]{article}
\usepackage{amsmath,amssymb,amsfonts,amsthm,url}
\usepackage[margin=1in]{geometry}

\usepackage{graphicx}
\usepackage{subcaption}

\theoremstyle{plain}
\newtheorem{theorem}{Theorem}[section]
\newtheorem{proposition}[theorem]{Proposition}
\newtheorem{corollary}[theorem]{Corollary}
\newtheorem{lemma}[theorem]{Lemma}

\theoremstyle{remark}
\newtheorem*{remark}{Remark}

\newcommand{\SC}{\widehat{K}}
\newcommand{\Z}{\mathbb{Z}}

\title{Settling the genus of the $n$-prism}
\author{Timothy Sun\\Department of Computer Science\\San Francisco State University}
\date{}

\begin{document}

\maketitle

\begin{abstract}
In a 1977 paper, Ringel conjectured a formula for the genus of the $n$-prism $K_n\times K_2$ and verified its correctness for about five-sixths of all values $n$. We complete this calculation by showing that, with the exception of $n = 9$, Ringel's conjecture is true for $n \equiv 5, 9\pmod{12}$. 
\end{abstract}

\section{Introduction}

Proving the Map Color Theorem~\cite{RingelYoungs,Ringel-MapColor} of Ringel, Youngs, \emph{et al.}\ amounts to demonstrating that the complete graph $K_n$ can be drawn without crossings in the surface of genus
$$\left\lceil\frac{(n-3)(n-4)}{12}\right\rceil,$$
matching a lower bound derived from the Euler polyhedral equation. As one of several possible followups to this classic result, one can ask about the genus $\gamma(G)$ of other related families of dense graphs. One such direction is the family of $n$-prism graphs $K_n \times K_2$, the graphs formed by taking two copies of $K_n$ and adding the edges of a perfect matching between the two graphs. Ringel~\cite{Ringel-nPrism} calculated the following lower bound on the genus of these graphs:

\begin{proposition}[Ringel~\cite{Ringel-nPrism}]
For $n \geq 3$, $$\gamma(K_n\times K_2) \geq \left\lceil\frac{(n-2)(n-3)}{6}\right\rceil.$$
\label{prop-lowerbound}
\end{proposition}

He went on to describe embeddings that matched this lower bound for all $n \not\equiv 5, 9 \pmod{12}$ using previously known embeddings, conjecturing that the lower bound is attained in these missing cases, except when $n = 5$. In this anomalous graph $K_5 \times K_2$, there are two disjoint copies of the nonplanar graph $K_5$, so its genus must be at least 2. We show, with computer assistance, that $n = 9$ is the only other counterexample and construct embeddings matching the lower bound for all the remaining graphs, yielding the following formula:

\begin{theorem}
For all $n \geq 2$, the genus of $K_n \times K_2$ is

\begin{equation*}
\gamma(K_n\times K_2) = 
\begin{cases} 
      \vspace{-13pt} \\
      \hspace{4pt}\displaystyle 
      \left\lceil\frac{(n-2)(n-3)}{6}\right\rceil+1 & 
                       \text{if~} n=5\text{~or~}9, \text{~and} 
      \\[16pt]
      \hspace{4pt}\displaystyle
      \left\lceil\frac{(n-2)(n-3)}{6}\right\rceil   & 
                       \text{otherwise.} 
      \\ \vspace{-12pt}
\end{cases}
\end{equation*}
\label{thm-main}
\end{theorem}

The missing cases are solved by constructing triangular embeddings of so-called ``split-complete'' graphs using the theory of current graphs. Current graphs were used in the original proof of the Map Color Theorem (see Ringel \cite{Ringel-MapColor}), and they are still the primary tool used to study embeddings of the complete graph: for example, a recent proof of the nonorientable Map Color Theorem by Korzhik \cite{Korzhik-Simple} uses current graphs to construct minimum genus embeddings of more than one complete graph at a time, and there has been an ongoing line of work on current graphs that produce face 2-colorable embeddings of the complete graphs (e.g., McCourt \cite{McCourt-Biembedding} and Archdeacon \cite{Archdeacon-Heffter}). Lawrencenko \emph{et al.}~\cite{Lawrencenko-4Genus} use current graphs to find quadrilateral embeddings of complete or near-complete graphs. Such embeddings are related to an analogue of the Map Color Theorem for embedded graphs with only even-length faces, which was proven in full by Liu \emph{et al.}~\cite{Liu-EvenMapColor}. 

Except for some of the simplest cases of the Map Color Theorem (see, e.g., the recursive constructions of Bonnington \emph{et al.}~\cite{Bonnington-Exponential} or the design-theoretic construction of Grannell \emph{et al.}~\cite{Grannell-SurfaceEmbeddings}), current graphs are the only known method for constructing minimum genus embeddings of the complete graphs $K_n$ for all sufficiently large $n$. 

Ringel~\cite{Ringel-nPrism} used embeddings of split-complete graphs to prove the genus formula of the prism graphs $K_n \times K_2$ when $n \equiv 8 \pmod{12}$, but an analogous construction for $n \equiv 5, 9 \pmod{12}$ using current graphs is not possible, since the vertex counts are of the opposite parity. We make use of two different approaches, both of which start from triangular embeddings of $K_{n+1}-K_3$.

\section{Background}

For background on topological graph theory, especially on the Map Color Theorem, see Gross and Tucker~\cite{GrossTucker} and Ringel~\cite{Ringel-MapColor}. 

The \emph{$n$-prism graph} is the graph $K_n \times K_2$, the Cartesian product of $K_n$ and $K_2$. More specifically, it is the graph with vertex set $\Z_n \times \Z_2$, where two vertices $(a_1,b_1)$ and $(a_2,b_2)$ are adjacent if and only if $a_1=a_2$ or $b_1=b_2$. We call the edges of the form $((a, 0),(a, 1))$ \emph{matching} edges. 

An (orientable) \emph{embedding} of a graph $G$ is an injective mapping $\phi\colon G \to S_g$ of the graph into the orientable surface of genus $g$, i.e., a sphere with $g$ handles. We restrict our attention to \emph{cellular} embeddings, those where the complement of the embedding $S_g \setminus \phi(G)$ is a disjoint union of disks. Each disk is referred to as a \emph{face} of the embedding. Embeddings of this kind are governed by the \emph{Euler polyhedral equation}
$$|V|-|E|+|F| = 2-2g,$$
where $V$, $E$, and $F$ denote the sets of vertices, edges, and faces, respectively. Each face can be thought of as a $k$-sided polygon, and a standard fact about cellular embeddings is that the total number of sides summed over all faces is $2|E|$ because each edge is incident with two sides.

Let $G$ be a graph, possibly with parallel edges or self-loops. We associate, with each edge $e \in E$, two opposing arcs $e^+$ and $e^-$, and we write $E^+$ to denote the set of all such arcs. A \emph{rotation} of a vertex $v$ is a cyclic permutation of all the arcs leaving $v$, and a \emph{rotation system} consists of a rotation for each vertex. When $G$ is simple, it is enough to give a cyclic permutation of each vertex's neighbors. Each cellular embedding gives rise to a rotation system by considering the clockwise ordering of the arcs in $E^+$ leaving a vertex with respect to some orientation of the surface. Likewise, every rotation system can be transformed into a cellular embedding by piecing together the embedded surface via face-tracing (see \S3.2.6 of Gross and Tucker~\cite{GrossTucker}). Each arc $e^\pm$ corresponds to a side of one face in the embedding. 

The genus of an embedding $\phi\colon G \to S_g$ simply refers to the genus $g$ of the embedded surface. The \emph{(minimum) genus} of a graph $G$, denoted by $\gamma(G)$, is the smallest integer $g$ for which there is an embedding in $S_g$. 

\subsection{The genus lower bound for prism graphs}

We begin by reproving Proposition~\ref{prop-lowerbound}, the lower bound on $\gamma(K_n \times K_2)$, since we are primarily interested in embeddings of genus exactly $(n-2)(n-3)/6$ (without rounding) that match this lower bound. Our argument is essentially the same as Ringel~\cite{Ringel-nPrism}, with some streamlined calculations. 

A standard approach to finding a lower bound for the genus of dense graphs is to use the girth to lower bound the length of any face. We note, however, that the $n$-prism graphs cannot have triangular embeddings, despite having girth 3, because the $n$ matching edges do not belong to any 3-cycle. Let $F_M$ denote the set of faces incident with at least one matching edge. For every face in $F_M$,
\begin{itemize}
\item there are at least two sides corresponding to matching edges, and
\item those instances are nonconsecutive.
\end{itemize}
The first property can be shown by observing that, since the matching edges form a cut, any closed walk passes through an even number of matching edges. Because there are $2n$ total sides that come from matching edges, this property implies that that $|F_M| \leq n$. 

The second property is the observation that the matching edges are vertex-disjoint. If we walk around each face boundary, each incidence of a matching edge is immediately followed by that of a non-matching edge, so there are at least an additional $2n$ sides in $F_M$ besides those from matching edges. The number of sides in $F_M$ is thus at least $4n$. Since there are $2\binom{n}{2}+n = n^2$ edges in the $n$-prism, the total number of sides remaining in $F\setminus F_M$ is at most $2n^2-4n$. Here, we apply the aforementioned girth argument to show that there can be at most $(2n^2-4n)/3$ faces outside of $F_M$. These two inequalities on $|F_M|$ and $|F\setminus F_M|$ can be merged into the upper bound $|F| \leq (2n^2-n)/3$. Substituting this inequality into the Euler polyhedral equation yields:
$$2 - 2g = |V|-|E|+|F| \leq 2n-n^2+(2n^2-n)/3,$$
which can be rearranged into
$$g \geq \frac{(n-2)(n-3)}{6}.$$

The upper bound on the number of faces is achieved exactly when each face in $F_M$ is quadrangular, and all other faces are triangular. For brevity, we call such an embedding of a prism graph \emph{snug}. 

\begin{corollary}
The genus of $K_n \times K_2$ is exactly $\frac{(n-2)(n-3)}{6}$ if and only if it has a snug embedding.
\label{cor-snug}
\end{corollary}

Another way of phrasing our main result is that there exist snug embeddings of $K_n \times K_2$ when $n \equiv 5 \text{~or~} 9 \pmod{12}$, and $n > 9$. 

\subsection{Current graphs}\label{sec-current}

Like in the proof of the Map Color Theorem, current graphs will be used to construct the desired genus embeddings. A \emph{current graph} $(G, \phi, \alpha)$ consists of an embedding of a graph $\phi\colon G \to S$ in an orientable surface $S$, along with an arc-labeling $\alpha\colon E^+ \to \Z_m$ such that $\alpha(e^+) = -\alpha(e^-)$ for all pairs of opposite arcs. In all of our drawings, black and white vertices represent clockwise and counterclockwise rotations, respectively. $\Z_m$ is known as the \emph{current group} and its elements are \emph{currents}. The \emph{excess} of a vertex $v$ is the sum of the currents over all arcs leaving $v$, and when the excess is 0, we say that the vertex satisfies Kirchhoff's current law. A vertex with nonzero excess is known as a \emph{vortex}.

The \emph{index} of a current graph refers to the number of faces in the embedding $\phi$. For each one of these faces, its boundary walk is called a \emph{circuit} and consists of a cyclic sequence of arcs $(e_1^\pm, e_2^\pm, \dotsc, e_i^\pm)$. The \emph{log} of the circuit is $(\alpha(e_1^\pm), \alpha(e_2^\pm), \dotsc, \alpha(e_i^\pm))$, the cyclic sequence of currents on those arcs. 

For the $n \equiv 9 \pmod{12}$ case, we use index 1 current graphs that satisfy these ``construction principles'':
\begin{itemize}
    \item [(A1)] Each vertex has degree 1 or 3.
    \item [(A2)] Each nonzero element of the current group $\Z_m$ appears in the log of the circuit exactly once.
    \item [(A3)] For each vertex of degree 3, Kirchhoff's current law is satisfied.
    \item [(A4)] The excess of each vertex of degree $1$ generates $\Z_m$. 
\end{itemize}
In other words, (A3) and (A4) imply that the vortices are exactly the degree 1 vertices. 

Index 3 current graphs are used for handling the $n \equiv 5\pmod{12}$ case, and they satisfy a similar set of properties:
\begin{itemize}
    \item [(B1)] Each vertex has degree 3.
    \item [(B2)] Each nonzero element of the current group $\Z_{3m}$ appears in the log of each circuit exactly once.
    \item [(B3)] For every vortex, each circuit is incident with that vortex, and its excess generates the index 3 subgroup $\langle 3 \rangle$ of $\Z_{3m}$, i.e., the subgroup generated by $3 \in \Z_{3m}$. 
    \item [(B4)] For each arc $e^+ \in E^+$, if circuit $[i]$ traverses $e^+$ and circuit $[j]$ traverses its opposite arc $e^-$, then $j-i \equiv \alpha(e^+) \pmod{3}$.
\end{itemize}

Associated with each current graph is a \emph{derived embedding} of a graph whose vertex set is $\Z_n$. In this paper, because all of our current graphs satisfy properties (A2) and (B2), the derived graph will always be the complete graph $K_n$ (see Gross and Tucker \cite{GrossTucker} for examples of non-complete derived graphs). Given an index $j$ current graph with circuits labeled $[0], [1], \dotsc, [j{-}1]$, the rotation at vertex $v \in \Z_n$ in the derived embedding is obtained by taking the log of circuit $[v \bmod{j}]$ and adding $v$ to each element. 

\begin{figure}[!t]
\centering
\includegraphics[scale=1]{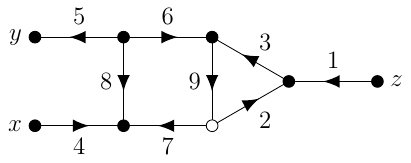}
\caption{An index 1 current graph.}
\label{fig-c9s1}
\end{figure}

When these construction principles are satisfied, the derived embedding has one Hamiltonian face for each vortex, and all other faces are triangular (see Ringel~\cite[\S2.3 and \S9.1]{Ringel-MapColor}). Each Hamiltonian face is then subdivided with a vertex, yielding a triangular embedding of $K_{n+\ell}-K_\ell$, where $\ell$ is the number of vortices. In our diagrams, each vortex is labeled with a distinct lowercase letter, and in the derived embedding, the corresponding subdivision vertex is also given that label. An example of an index 1 current graph satisfying these properties is given in Figure~\ref{fig-c9s1}, where there are three vortices and the current group is $\mathbb{Z}_{19}$. The log of its only circuit is
$$\begin{array}{rrrrrrrrrrrrrrrrrrrrrrrrrrrrrrrrrrrrrrrrrrrrr}
\lbrack0\rbrack. & 15 & x & 4 & 11 & 5 & y & 14 & 6 & 16 & 18 & z & 1 & 17 & 10 & 13 & 8 & 12 & 2 & 3 & 9 & 7, \\ 
\end{array}$$
where we follow the usual convention of including vortex labels in the log where they appear along the circuit. The derived embedding, after subdivision, is a triangular embedding of $K_{22}-K_3$, from which we will later obtain a genus embedding of $K_{21}\times K_2$.

\section{Ringel's reflection}

We say that a subset of the faces of an embedded graph is a \emph{facial cover} if every vertex is incident with at least one of the faces. All of Ringel's constructions~\cite{Ringel-nPrism} can be summarized in the following way:

\begin{lemma}
Let $\phi\colon K_n \to S_g$ be an embedding in the surface of genus $g$, and let $F_c \subseteq F$ be a facial cover of size $h$. Then there exists an embedding of $K_n \times K_2$ in the surface of genus $2g+h-1$. 
\label{lem-mainconstruction}
\end{lemma}

\begin{proof}
We start by making a mirror image $\phi'\colon K_n \to S_g$ of the embedding $\phi$. For each cell (i.e., face, vertex, etc.) $c$, we label its corresponding mirror image $c'$. The facial cover $F_c$ induces a mirror image $F_c'$ in $\phi'$. For each face $f \in F_c$ and its mirror image $f' \in F_c'$, we puncture their interiors and connect the resulting boundaries with a tube. Along this tube we may add the edge $(v, v')$ for every vertex $v$ incident with $f$. The first tube forms the connected sum of the two surfaces, which has genus $2g$, and each of the additional $h-1$ tubes becomes a handle that increases the genus by 1. 
\end{proof}

\begin{figure}[!t]
\centering
\includegraphics[scale=0.93]{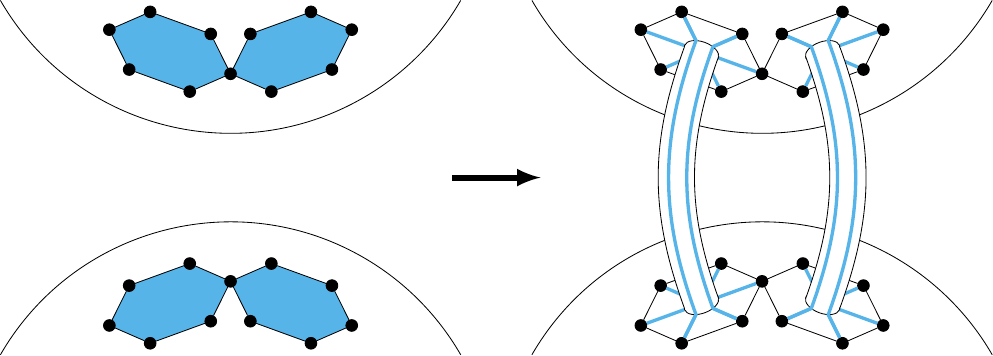}
\caption{Adding tubes between two mirror images to create a prism.}
\label{fig-mainconst}
\end{figure}

The construction is illustrated in Figure~\ref{fig-mainconst} with a facial cover of size 2. We note that for the vertex and its mirror image incident with both faces in their respective facial covers, only one edge is added, resulting in a larger, non-quadrilateral face along the other tube. Broadly speaking, Lemma~\ref{lem-mainconstruction} is most economical, in terms of the genus of resulting embedding, when the faces of the facial cover overlap as little as possible.

For many of the prism graphs $K_n \times K_2$, it suffices to use a genus embedding of $K_{n+1}$.

\begin{corollary}[Ringel~\cite{Ringel-nPrism}]
There exists an embedding of $K_n \times K_2$ in the surface of genus $$2\left\lceil\frac{(n-2)(n-3)}{12}\right\rceil.$$
\label{cor-singlehandle}
\end{corollary}
\begin{proof}
Starting from a genus embedding of $K_{n+1}$, delete a vertex and let the facial cover consist of just the resulting face. 
\end{proof}

One can verify that for $n \equiv 2, 3, 6, 7, 10, 11 \pmod{12}$, the genus of this embedding matches that in Theorem~\ref{thm-main}. For the other residues, the genus of these embeddings are off by one from Ringel's conjecture, necessitating alternative constructions involving larger facial covers. When $n \equiv 1, 4 \pmod{12}$, facial covers with overlapping vertices can still be used:

\begin{corollary}[Ringel~\cite{Ringel-nPrism}]
For $n \equiv 1, 4\pmod{12}$, there exists an embedding of $K_n \times K_2$ in the surface of genus
$$\left\lceil\frac{(n-2)(n-3)}{6}\right\rceil.$$
\end{corollary}

\begin{proof}
We begin with triangular embeddings of $K_{n+1}-K_2$, which were constructed in, e.g., Youngs~\cite{Youngs-3569} and Jungerman~\cite{Jungerman-KnK2}. Let $(a,b)$ be the missing edge. By deleting, say, vertex $a$, we obtain a large face incident with all vertices besides $b$. Then, a facial cover of size 2 can be formed by combining this face with any face incident with $b$. If $g$ is the genus of $K_{n+1}-K_2$, then the genus of $K_{n+1}$ must be $g+1$, since $K_{n+1}-K_2$ has a triangular embedding. Thus, the embedding of Corollary \ref{cor-singlehandle} has genus $2(g+1)$, but this new embedding has genus $2g+1$ and hence is optimal. 
\end{proof}

For the remaining residues $n \equiv 0, 5, 8, 9$, a snug embedding is required. We first characterize when Lemma~\ref{lem-mainconstruction} can be used to generate such embeddings. Let us call a facial cover $F_c$ a \emph{cotriangular patchwork} if each vertex has exactly one incidence with the faces in $F_c$ and every face in $F\setminus F_c$ is triangular. 

\begin{proposition}
Let $\phi\colon K_n \to S_g$ be an embedding, and let $F_c \subseteq F$ be a cotriangular patchwork. Then applying Lemma~\ref{lem-mainconstruction} to $\phi$ and $F_c$ yields a snug embedding of $K_n \times K_2$.
\label{prop-patchwork}
\end{proposition}
\begin{proof}
Each matching edge is inserted in exactly one tube, so each tube consists of all quadrangular faces. All the other faces are triangular faces that were already present in the original embedding.
\end{proof}

\begin{remark}
The approach described here is a slight variation on the White-Pisanski method \cite{White-RepeatedCartesian, Pisanski-Cartesian} (also see Gross and Tucker \cite[\S3.5.4]{GrossTucker}) for constructing quadrilateral embeddings of Cartesian products. The difference is that we require all other faces outside of the patchwork to be triangular instead of quadrilateral.
\end{remark}

\begin{corollary}[Ringel~\cite{Ringel-nPrism}]
There exists a snug embedding of $K_{12s} \times K_2$, $s \geq 1$.
\end{corollary}
\begin{proof}
The triangular embeddings of $K_{12s}$ of Terry \emph{et al.}~\cite{Terry-Case0} have cotriangular patchworks consisting of $4s$ triangular faces.
\end{proof}

Let a graph $\SC_n$ be called a \emph{split-complete graph} if one can label the vertices with numbers $1, 2, \dotsc, n-1$ and letters $u$ and $v$ such that the numbered vertices are all pairwise adjacent, and the neighbors of $u$ and $v$ partition the numbered vertices. Ringel used a family of triangular embeddings of split-complete graphs $\SC_{12s+9}$ (see Ringel \cite[\S6.9]{Ringel-MapColor}) to find genus embeddings of $K_{12s+8}\times K_2$. A triangular embedding of a split-complete graph $\SC_{n+1}$ is equivalent to an embedding of the complete graph $K_n$ with a cotriangular patchwork by deleting $u$ and $v$ from the embedding of $\SC_{n+1}$.

\begin{corollary}[Ringel~\cite{Ringel-nPrism}]
There exists a snug embedding of $K_{12s+8} \times K_2$, $s \geq 0$.
\end{corollary}

We note that another solution for $n \equiv 0 \pmod{12}$ is possible using the embeddings of split-complete graphs $\SC_{12s+1}$ constructed by Jungerman \cite{Jungerman-OrientableCascades} and the author \cite{Sun-Index2}. The common trait shared between the aforementioned families of split-complete graphs is that the special vertices $u$ and $v$ are each incident with exactly half of the numbered vertices. However, for $n \equiv 5, 9\pmod{12}$, there are an odd number of such vertices, so alternative approaches are needed. On a more technical level, those families of current graphs used a vortex whose excess generates the index 2 subgroup of an even-order cyclic group.

\section{A counterexample}

The lone exception to Ringel's conjecture is the case of $n=9$. The following fact allows us to restrict our attention to just embeddings of $K_9$:

\begin{proposition}
Suppose there exists a snug embedding of $K_n \times K_2$. Then the embedding restricted to either copy of $K_n$ has a cotriangular patchwork. 
\label{prop-slice}
\end{proposition}
\begin{proof}
Each of the matching edges lie on two quadrilateral faces, so their duals form a disjoint union of cycles. Cutting the surface along these dual edges splits the surface into two punctured surfaces, each embedded with a copy of $K_n$. The punctured faces of each surface form a cotriangular patchwork.
\end{proof}

\begin{proposition}
The genus of $K_9 \times K_2$ is $8$. 
\label{prop-n9}
\end{proposition}
\begin{proof}
Corollary~\ref{cor-singlehandle} produces an embedding of $K_9 \times K_2$ of genus 8. If there were an embedding of $K_9 \times K_2$ of genus 7, by Corollary~\ref{cor-snug} and Proposition~\ref{prop-slice}, its restriction to each copy of $K_9$ would have a cotriangular patchwork. Furthermore, at least one of those restrictions must be a genus embedding in $S_3$, or else the genus of the entire embedding would be at least 8. We focus on this particular genus embedding $K_9 \to S_3$. 

Any nontriangular face must be part of the cotriangular patchwork. The Euler polyhedral equation forces the set of nontriangular faces to be of one of three varieties:
\begin{enumerate}
\item[(a)] three 4-sided faces,
\item[(b)] one 5-sided face and one 4-sided face, or
\item[(c)] one 6-sided face.
\end{enumerate}
In case (a), those faces cannot be a cotriangular patchwork because the pigeonhole principle guarantees there to be a vertex with more than one incidence. For case (b), these two faces are the patchwork, and for case (c), the patchwork is the 6-sided face joined by one of the triangular faces. A computer search shows that no embedding of $K_9$ has either type of patchwork. 
\end{proof}
The above result can also be verified using the exhaustive enumerations of triangulations of $S_3$ by Sulanke and Lutz~\cite{SulankeLutz-Isomorphism}.\footnote{Their lists are available at \url{http://page.math.tu-berlin.de/~lutz/stellar/surfaces.html}}

\section{The current graph constructions}

We follow previous work (see Ringel~\cite{Ringel-MapColor}) by using current graphs to construct triangular embeddings of $K_{n+1}-K_3$, $n \equiv 5, 9 \pmod{12}$. Unlike in the proof of the Map Color Theorem, we enforce additional constraints in order to transform them into triangular embeddings of split-complete graphs.

There are several known families of index 3 current graphs that generate triangular embeddings of $K_{12s+6}-K_3$ (see, e.g., \cite{Sun-Minimum} for the simplest known constructions). One local operation for modifying a triangular embedding is the \emph{edge flip}, which removes an edge and adds in the other diagonal of the resulting quadrilateral face. When vortices are adjacent in the current graph, edge flips on the derived embedding can create triangular embeddings of other graphs with the same number of edges. One of those graphs, $K_{12s+6}-K_{1,3}$, is particularly useful for our purposes because all the missing edges involve the same vertex $v$. If the non-neighbors of $v$ belong to the same triangular face, then we would obtain a split-complete graph $\SC_{12s+6}$ by subdividing that triangle, or equivalently, deleting vertex $v$ yields a cotriangular patchwork. Triangular embeddings of $K_{12s+6}-K_{1,3}$, $s \geq 3$ were first found by Gross~\cite{Gross-Case6}, based on an earlier construction of Youngs~\cite{Youngs-3569}. However, Gross's method never results in an embedding satisfying our desired requirements on the non-neighbors. We require both a different family of current graphs and a different sequence of edge flips.

The aforementioned edge flip approach cannot be applied to the simple families of index 1 current graphs for $K_{12s+10}-K_3$ (see \cite{Youngs-1710} or \cite[\S2.3]{Ringel-MapColor}) since the degree 1 vortices cannot be adjacent. We circumvent this problem by adding the missing three edges using a handle and later ``subtracting'' a different handle to arrive at a split-complete graph. 

\subsection{$n \equiv 5 \pmod{12}$}

\begin{figure}[!t]
\centering
\includegraphics[scale=0.93]{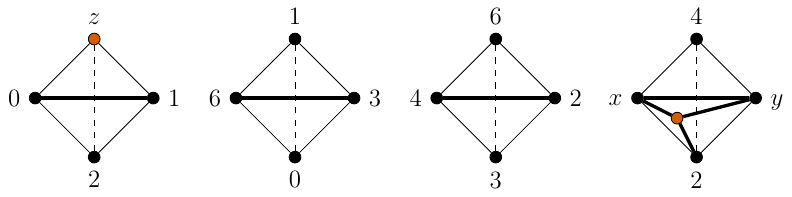}
\caption{Flipping edges and subdividing one face to make a split-complete graph.}
\label{fig-edgeflip6}
\end{figure}

\begin{theorem}
There exists a triangular embedding of a split-complete graph $\SC_{12s+6}$, for all $s \geq 1$.
\label{thm-c5}
\end{theorem}
\begin{proof}
For $s = 1$, a triangular embedding of a split-complete graph $\SC_{18}$ is given in Mayer \cite{Mayer-Orientables}.

\begin{figure}[!t]
\centering
\includegraphics[scale=0.93]{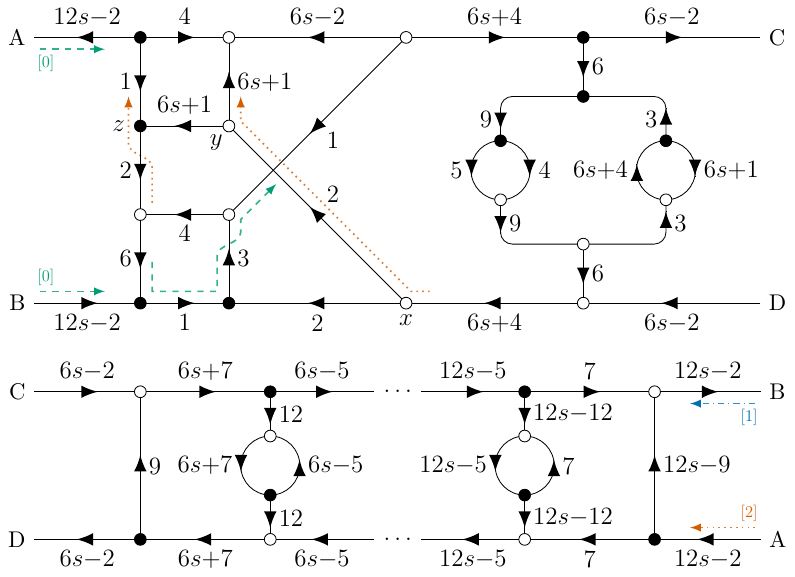}
\caption{A family of index 3 current graphs with current group $\Z_{12s+3}$.}
\label{fig-current6}
\end{figure}

The family of current graphs in Figure~\ref{fig-current6} generates triangular embeddings of $K_{12s+6}-K_3$ for $s \geq 2$. The construction is similar to previous index 3 constructions by the author~\cite{Sun-Minimum}, except that the repeating pattern in the second half of the picture is derived from Youngs's solution to ``Case 3'' of the Map Color Theorem, instead of ``Case 5'' (see Youngs~\cite{Youngs-3569} or Ringel~\cite[\S9.2]{Ringel-MapColor}). The partial circuits in Figure \ref{fig-current6} indicate that the logs of the circuits are of the form 
$$\begin{array}{rrrrrrrrrrrrrrrrrrrrrrrrrrr}
\lbrack0\rbrack. & \dots \,\, 6 & 1 & 3 & -1 \,\, \dots \\ 
\end{array}$$
$$\begin{array}{rrrrrrrrrrrrrrrrrrrrrrrrrrr}
\lbrack2\rbrack. & \dots \,\, x & 2 & y \,\, \dots \,\, -2 & z & -1 \,\, \dots
\end{array}$$
Thus, in the derived embedding, the rotations of vertices $0$, $2$, and $3$ are of the form
$$\begin{array}{rrrrrrrrrrrrrrrrrrrrrrrrrrr}
0. & \dots & 6 & 1 & 3 & \dots & \\
2. & \dots & 0 & z & 1 & \dots & x & 4 & y & \dots \\
3. & \dots & 4 & 6 & 2 & \dots &
\end{array}$$
Figure~\ref{fig-edgeflip6} depicts a sequence of edge flips at these locations in the derived embedding, where vertical dashed edges are replaced with thick horizontal edges. The resulting embedding is of $K_{12s{+}6}-K_{1,3}$, where the missing edges are $(z, 2)$, $(z,x)$, and $(z,y)$. Furthermore, the last edge flip creates the triangular face $\lbrack 2, y, x\rbrack$, so subdividing that face with a new vertex yields a triangular embedding of a split-complete graph $\SC_{12s{+}6}$.
\end{proof}

\begin{corollary}
There exists a snug embedding of the prism graph $K_{12s+5}\times K_2$, for all $s \geq 1$.
\label{cor-c5}
\end{corollary}

\subsection{$n \equiv 9 \pmod{12}$}

To transform a triangular embedding of a split-complete graph $\SC_{n+1}$ into a genus embedding of $K_{n+1}$, the edge $(u,v)$ is first inserted arbitrarily into the rotation system and then contracted to create a single vertex incident with all others. The two triangular faces merged by this operation turn into a hexagonal face that has two incidences with the identified vertex. We rely on the fact that this process can be reversed:

\begin{lemma}
Suppose there exists an embedding $\phi\colon K_{n+1} \to S_g$ with exactly one nontriangular face, and that face is hexagonal and incident with some vertex twice. Then there exists a triangular embedding of a split-complete graph $\SC_{n+1}$ in the surface $S_{g-1}$.
\label{lem-split}
\end{lemma}
\begin{proof}
Since the embedding is orientable, the hexagonal face is of the form $\lbrack w, a, b, w, c, d\rbrack$ where $a, b, c, d$ are all distinct (see Lemma 2.5 of Sun~\cite{Sun-FaceDist}). Suppose the rotation at $w$ is of the form
$$\begin{array}{rrrrrrrrrrrrrrrrrrrrrrrrrrr}
w. & d & a & p_1 \,\, \dots \,\, p_i & b & c & q_1 \,\, \dots \,\, q_j,
\end{array}$$
where $p_1, \dotsc, p_i, q_1, \dotsc, q_j$ are the remaining vertices. Following Figure~\ref{fig-splitting}, we split vertex $w$ by replacing it with two adjacent vertices $u$ and $v$ with rotations
$$\begin{array}{rrrrrrrrrrrrrrrrrrrrrrrrrrr}
u. & v & a & p_1 \,\, \dots \,\, p_i & b
\end{array}$$
and
$$\begin{array}{rrrrrrrrrrrrrrrrrrrrrrrrrrr}
v. & u & c & q_1 \,\, \dots \,\, q_j & d.
\end{array}$$

\begin{figure}[!t]
\centering
\includegraphics[scale=0.93]{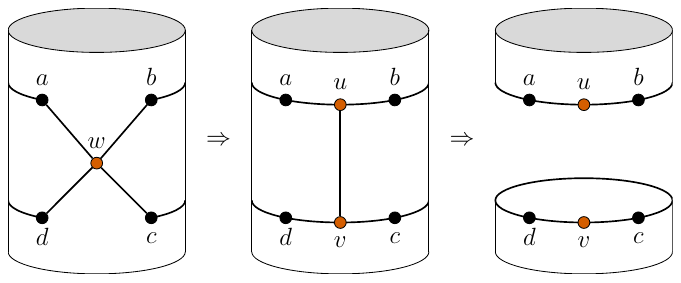}
\caption{Splitting a doubly-incident vertex to remove a handle.}
\label{fig-splitting}
\end{figure}

The original 6-sided face is now 8-sided. Since the edge $(u,v)$ is incident with the 8-sided face twice, its deletion causes the face to split into two triangles $\lbrack u, a, b\rbrack$ and $\lbrack v, c, d \rbrack$, thereby decreasing the genus by one. Since the original embedding was of a complete graph and had its one nontriangular face broken into triangles, the final embedding is triangular and of a split-complete graph.
\end{proof}

In previous work by the author~\cite{Sun-FaceDist}, genus embeddings of $K_{n+1}$ with at most one nontriangular face were found for all $n \neq 7$. We describe a similar construction for the $n \equiv 9 \pmod{12}$ case, with additional constraints to ensure that the hexagonal face has a repeated vertex. 

\begin{figure}[!t]
\centering
    \begin{subfigure}[b]{\textwidth}
        \centering
        \includegraphics[scale=0.93]{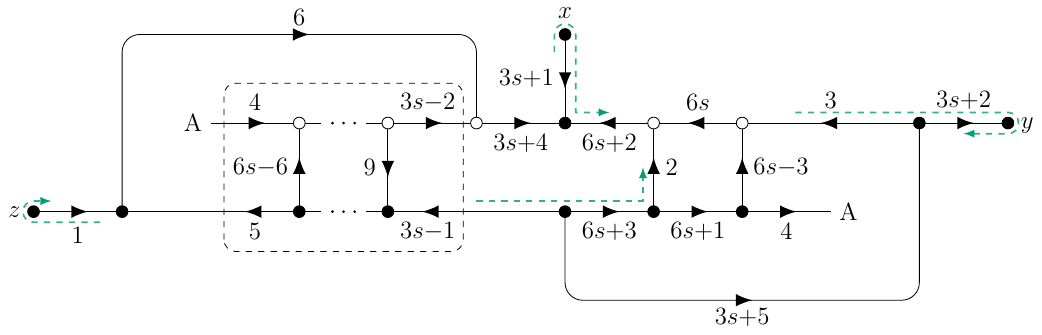}
        \caption{}
    \end{subfigure}
    \begin{subfigure}[b]{\textwidth}
        \centering
        \includegraphics[scale=0.93]{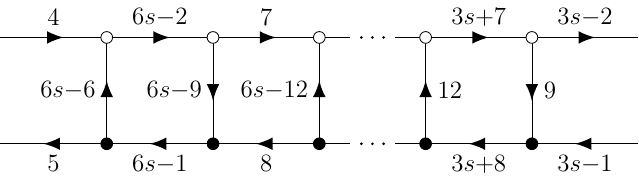}
        \caption{}
    \end{subfigure}
\caption{A family of index 1 current graphs with current group $\Z_{12s+7}$ (a) and an expanded view of the ladder of varying length (b).}
\label{fig-current10}
\end{figure}

\begin{figure}[!t]
\centering
    \begin{subfigure}[b]{0.7\textwidth}
        \centering
        \includegraphics[scale=0.85]{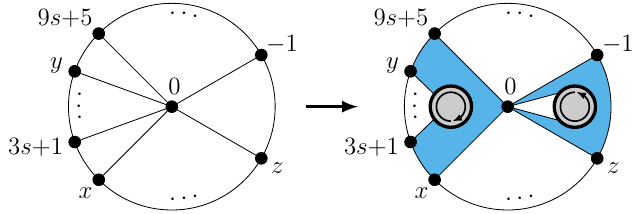}
        \caption{}
    \end{subfigure}
    \begin{subfigure}[b]{0.29\textwidth}
        \centering
        \includegraphics[scale=0.85]{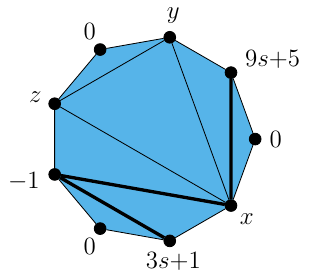}
        \caption{}
    \end{subfigure}
    \begin{subfigure}[b]{0.7\textwidth}
        \centering
        \includegraphics[scale=0.85]{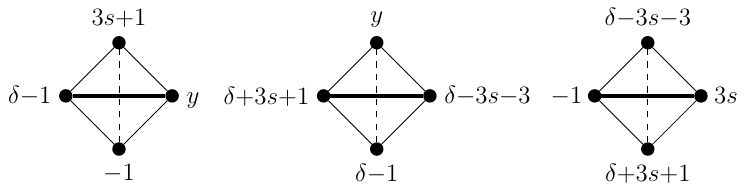}
        \caption{}
    \end{subfigure}
    \begin{subfigure}[b]{0.29\textwidth}
        \centering
        \includegraphics[scale=0.85]{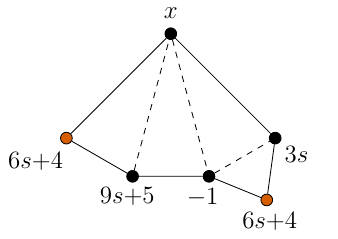}
        \caption{}
    \end{subfigure}
\caption{In the derived embedding, adding a handle near vertex 0 (a) creates a 9-sided face where we can add edges (b), and a few additional edge flips (c) produces the desired non-simple hexagonal face (d).}
\label{fig-edgeflip10}
\end{figure}

\begin{theorem}
There exists a triangular embedding of a split-complete graph $\SC_{12s+10}$, for all $s \geq 1$.
\label{thm-c9}
\end{theorem}
\begin{proof}
Consider the current graphs in Figures~\ref{fig-c9s1} and \ref{fig-current10}. These represent the $s = 1$ and $s \geq 2$ cases, respectively. We can follow the log in Section \ref{sec-current} or, like in the previous construction, the partial circuits in Figure \ref{fig-current10} to see that the log of the single circuit in each of the current graphs is of the form
$$\begin{array}{rrrrrrrrrrrrrrrrrrrrrrrrrrr}
\lbrack0\rbrack. & x & 3s{+}1 & 6s{+}5 \,\, \dots \,\, \delta & 3s{+}2 & y & 9s{+}5 \,\, \dots \,\, {-}1 & z \,\, \dots \,\, 9s{+}5{-}\delta & 6s{+}3 & -1{-}\delta \,\, \dots
\end{array}$$ 
where $\delta = 11$ for $s = 1$ and $\delta = -3$ for $s \geq 2$. The value $\delta$ will determine slightly different sequences of edge flips between these two cases. In the derived embedding, the log is equal to the rotation at vertex 0 and induces the following partial rotations:
$$\begin{array}{rcccccccccccccccccccccccc}
-1. & \dots & \delta-1 & 3s{+}1 & y & \dots & x & 3s & 6s{+}4 \\
\delta{-}1. & \dots & \delta{+}3s{+}1 & y & \delta{-}3s{-}3 & \dots \\
\delta{+}3s{+}1. & \dots & -1 & \delta{-}3s{-}3 & 3s & \dots \\
x. & \dots & 3s & -1 & 9s{+}5 & 6s{+}4 & \dots
\end{array}$$
By rearranging the rotation at vertex 0, as shown in Figure~\ref{fig-edgeflip10}(a), we draw in the missing edges $(x,y)$, $(x,z)$, and $(y,z)$ and some other duplicate edges as in Figure~\ref{fig-edgeflip10}(b). After flipping the edges shown in Figure~\ref{fig-edgeflip10}(c), we can delete one copy of each duplicate edge to obtain the 6-sided face in Figure~\ref{fig-edgeflip10}(d) with the vertex $6s{+}4$ repeated. This vertex can be split using Lemma~\ref{lem-split} to create the desired embedding of a split-complete graph $\SC_{12s+10}$, for all $s \geq 1$. 
\end{proof}

\begin{corollary}
There exists a snug embedding of the prism graph $K_{12s+9}\times K_2$, for all $s \geq 1$.
\label{cor-c9}
\end{corollary}

\bibliographystyle{alpha}
\bibliography{biblio}

\end{document}